%

%
\magnification=\magstep1
\input amstex
\documentstyle{amsppt}

\def\c{\cite}
\def\l{\langle}
\def\r{\rangle}
\def\st{such that }
\def\iff{if and only if }
\def\O{OTTER }
\leftheadtext{OTTER Experiments in a System of Combinatory Logic}
\NoBlackBoxes

\topmatter
\title
OTTER Experiments in a System of Combinatory Logic$^*$
\endtitle
\author
Thomas Jech
\endauthor
\affil
Mathematics Department \\
The Pennsylvania State University  \\
University Park, PA 16802
\endaffil
\address
Mathematics Department,
The Pennsylvania State University,
University Park, PA 16802
\endaddress
\email jech\@math.psu.edu \endemail
\abstract
This paper describes some experiments involving the automated
theorem-proving program \O in the system TRC of illative 
combinatory logic.  We show how \O can be steered to find a 
contradiction in an inconsistent variant of TRC, and present
some experimentally discovered identities in TRC.
\endabstract
\endtopmatter

\document
\footnote"{}"{$^*$This 
work was supported by the Applied Mathematical Sciences program of the 
U. S. Department of Energy, under the Summer 1994 Faculty Research
Participation Program at Argonne National Laboratory.}

\specialhead
1. \ Introduction.
\endspecialhead

\O \c5 is a resolution/paramodulation theorem-proving program for
first-order logic with equality.  It has been used successfully in
several areas of logic and algebra \c8, \c9, \c6, \c7, \c4.

In this paper we describe our experiments with \O in the system
TRC of illative combinatory logic \c3.  The system TRC has been formulated 
by M. R. Holmes who proved that it is equiconsistent with Quine's {\it New 
Foundations.}

{\it New Foundations} (NF) was introduced by Quine and extensively
studied by Rosser, Jensen and others. NF is an attempt to axiomatize
mathematics in a way different from the traditional Zermelo-Fraenkel
axioms. Specker showed that NF refutes the Axiom of Choice, but the
consistency of NF is an open question.

Combinatory logic axiomatizes the intuition that everything is a function
{\it (combinator).} The basic symbol in the language is a binary
function $a$ and $a(x,y)$ (usually abbreviated to $x(y)$ or just
$xy$) denotes the result of applying $x$ to $y.$

Holmes gave a translation of NF into a combinatory logic system called TRC.
\O is particularly well suited for experimentation in TRC: \ TRC
has a small number of axioms, most of them equations, and apart from
its equiconsistency with NF, not much is known about the theory (we
wish to point out that no model of TRC is known).  This work (as well
as the ongoing experimentation by my student W. Wood) is a modest
beginning in a systematic investigation of TRC using automated 
reasoning techniques.

The theory TRC is closely related to combinatory logic \c2, \c1.  When
TRC is slightly modified one obtains an inconsistent theory.  This is
easily proved by elementary methods of lambda calculus (and well known
to the author of TRC).  Here we describe how a contradiction in the
inconsistent theory can be found directly, using \O.  

While experimenting with \O, we have come across a number of
interesting consequences of the axioms of TRC, and we present
here a few, with proofs provided by \O.

The work on this paper was done while the author visited Argonne
National Laboratory in May 1994, under the Faculty Research
Participation Program.  I am greatly indebted to William McCune and
Larry Wos for their patient explanations of the workings of \O, and
for many valuable discussions we had on automated reasoning
techniques.

\specialhead
The Theory TRC
\endspecialhead

M. R. Holmes introduced in \c3 the system TRC (for `type-respecting 
combinators) and proved that TRC is equiconsistent with Quine's New
Foundations.  The theory TRC is a system of combinatory logic: \ the
objects of TRC are {\it combinators}.  For two combinators $x$ and
$y$ we use concatenation to denote application
$$
	xy=x(y)
$$
and use the convention that the operation $xy$ associates left:
$$
	xyz=(xy)z.
$$

The theory TRC has four constant combinators Abst, Eq, $p_1$ and $p_2$,
the operation (function) of application, a two-place function pair$(x,y)$ 
(written as $\l x,y\r$), and a one-place function $K$.  The axioms of
TRC are as follows:
\roster
\item"{I. \ }"  $K(x)y=x$ \qquad\qquad  (constant functions)
\item"{II. \ }" $p_i\l x_1,x_2\r=x_i\quad$ for $i=1,2$ \qquad (projections)
\item"{III. \ }" $\l p_1 x,p_2 x\r=x$ \qquad\qquad (pairing)
\item"{IV. \ }" $\l f,g\r x=\l fx,gx\r$ \qquad\qquad (pairwise application)
\item"{V. \ }"  Abst $x\,y\,z=x\,K(z)(y\,z)$ \qquad\qquad (abstraction)
\item"{VI. \ }" Eq $\l x,y\r=p_1$ if $x=y$; Eq $\l x,y\r=p_2$  if $x\neq y$
\quad(equality)
\item"{VII. \ }" If for all $x$, $fx=gx$, then $f=g$ \qquad\quad(extensionality) 
\item"{VIII. \ }" $p_1\neq p_2$
\endroster

\noindent In \c3 Holmes proves that NF can be interpreted
in TRC, and conversely, TRC can be interpreted in NF.  As NF is as
yet not known to be either consistent or inconsistent, the same
applies to TRC.

\specialhead
An Inconsistent Variant of TRC
\endspecialhead

It is important that $K$ is a function symbol, not a combinator,
even though its definition is formally the same as that of the 
combinator $K$ in combinatory logic:   
\roster
\item"{I$^*$. \ }" $K\,x\,y=x.$
\endroster

If we replace the function symbol $K$ by a combinator $K$ and replace
Axiom I by I$^*$ and Axiom V by
\roster
\item"{V$^*$. \ }" $\text{Abst }x\,y\,z=x(Kz)(y\,z),$
\endroster
then the resulting system TRC$^*$ is inconsistent.  This can be proved 
by employing the simple but powerful feature of combinatory logic,
the {\it abstraction property}, along with a version of Russell's
paradox.

The actual explicit contradictory term might be quite complicated, and
a natural question arises whether an automated theorem prover can find
a contradiction, without being prompted by a human.

We have conducted a large number of experiments with \O and 
eventually succeeded to have \O find a contradiction, but not 
without a few hints.

\specialhead
The Strategy
\endspecialhead

First we restate the axioms of TRC$^*$ into clauses.  Even though \O
does accept first order formulas as input, it translates them into 
clauses anyway.  Besides, an inspection of the clauses that form
the input should help decide which options to select.
\bigskip
\eightpoint
{\tt\obeylines\obeyspaces
  
set(knuth\_bendix).
set(ur\_res).
set(unit\_deletion).
set(bird\_print).
  
assign(max\_mem,64000).
assign(max\_weight,40).                              
assign(pick\_given\_ratio,6).
  
list(sos).
x = x.
a(a(k,x),y) = x.         
a(p1,pair(x,y)) = x.     
a(p2,pair(x,y)) = y.          
pair(a(p1,x),a(p2,x)) = x.         
a(pair(x,y),z) = pair(a(x,z),a(y,z)).         
a(a(a(abst,x),y),z) = a(a(x,a(k,z)),a(y,z)).  
a(eq,pair(x,x)) = p1.             
x = y | a(eq,pair(x,y)) = p2.  
x = y | a(x,n(x,y)) != a(y,n(x,y)).         
p1 != p2.
end\_of\_list.

}
\tenpoint
\bigskip

Note that to state the axiom of extensionality one has to introduce a
Skolem function $n(x,y)$ $(n(x,y)=$ some $z$ such that $xz\neq yz$.
Also, we write  $a(x,y)$ for $xy$ (but the bird-print option makes it
possible to output $a$ as concatenation).  

Next, we decide which
clauses to  put in the set of support.  As we are looking for a
contradiction and do not intend to concentrate on any particular
clause, we put all clauses into {\tt sos}.

Most clauses in the axiom set are units, and in fact equations.  For
that reason we select the Knuth-Bendix completion procedure \c5.  This
procedure transforms a set of equalities into a set of rewrite rules.
For a suitable inference rule to accompany {\tt knuth\_bendix}, we follow
McCune's advice to use the unit-resulting resolution (UR-resolution),
along with unit deletion.  As our SPARCstation has sufficient memory,
we set the maximum to 64MB.  As we expect the contradiction to be
quite complicated, we set the maximum weight = 40.

\specialhead
Search for a Contradiction
\endspecialhead

With the input file and options set as described above, \O exceeded
the allocated memory without finding a proof.  It generated well over
100,000 clauses, mostly equations, but did not discover a
contradiction.

Our first attempt to enhance its power was to expand the input by
adding ``interesting'' clauses proved by \O, and introduce names
for ``important'' terms.  Without a clear direction, this approach
failed.  With the language being enlarged, the number of ``useless''
tautologies grew rapidly and \O reached the memory limit more
easily.  After many unsuccessful runs we decided to abandon this
strategy.  (This strategy however turned out to be a good exploratory
tool: \ \O generated a number of clauses that we later verified to
be theorems of TRC).

As it became clear that this ``formalist'' approach would not work we
had to decide how to steer \O in the right direction.  As a first
gentle push, we directed \O to search for a ``diagonal''
combinator.  A contradiction would no doubt involve some form of
diagonalization, so we asked \O to find a combinator $F$ with the
property $F\,x=x\,x$.

When this attempt was unsuccessful, we tried several variants,
eventually finding a combinator $F=\text{Abst Abst }K$ \st
$$
	F\,x\,y=x\,x.
$$
Below is a run of \O that found this $F$.  Since looking for a term
satisfying a single equation should not require substitutions into or
from nonatomic clauses, we set the flags {\tt para\_from\_units\_only} and
{\tt para\_into\_units\_only}.  This focused the search sufficiently and \O
found an answer in 378 sec. (after 22461 clauses).

\bigskip
\eightpoint

{\tt\obeylines\obeyspaces

set(knuth\_bendix).
set(ur\_res).
set(unit\_deletion).
set(para\_from\_units\_only).
set(para\_into\_units\_only).
set(bird\_print).
assign(max\_mem,64000).
assign(max\_weight,40).
assign(pick\_given\_ratio,6).

list(usable).
0 [] x=x.
end\_of\_list.

list(sos).
0 [] k x y=x.
0 [] p1 pair(x,y)=x.
0 [] p2 pair(x,y)=y.
0 [] pair(p1 x,p2 x)=x.
0 [] pair(x,y) z=pair(x z,y z).
0 [] abst x y z=x (k z) (y z).
0 [] eq pair(x,x)=p1.
0 [] x=y|eq pair(x,y) =p2.
0 [] x=y|x n(x,y) !=y n(x,y).
0 [] p1!=p2.
0 [] y b(y) c(y)!=b(y) b(y) |\$ans(y).
end\_of\_list.

----> UNIT CONFLICT at 378.27 sec ----> 22461 [binary,22460.1,1.1] 

\$ans(abst abst k).

---------------- PROOF ----------------

1 [] x=x.
3,2 [] k x y=x.
12 [] abst x y z=x (k z) (y z).
18 [] x b(x) c(x)!=b(x) b(x) |\$ans(x).
21 [] x (k y) (z y)=abst x z y.
104 [para\_from,12.1.1,18.1.1.1] 
 x (k b(abst x y)) (y b(abst x y)) c(abst x y)!=b(abst x y) b(abst x y) |
 \$ans(abst x y).
269 [para\_into,21.1.1.2,2.1.1] abst x (k y) z=x (k z) y.
22460 [para\_from,269.1.1,104.1.1,demod,3]
 b(abst abst k) b(abst abst k)!=b(abst abst k) b(abst abst k) |\$ans(abst abst k).
22461 [binary,22460.1,1.1] \$ans(abst abst k).

------------ end of proof -------------

}
\tenpoint
\bigskip

An inspection of 
the proof shows that only Axioms I$^*$ and V$^*$ are used.  When we
deleted the irrelevant information and let \O concentrate on $a,$
$Abst$ and $K$ only, it found the answer in less than a second.  This
illustrates how essential it is to choose the input to contain only 
as many assumptions as necessary for the proof.

\bigskip
\eightpoint

{\tt\obeylines\obeyspaces

set(knuth\_bendix).
set(ur\_res).
set(unit\_deletion).
set(para\_from\_units\_only).
set(para\_into\_units\_only).
set(bird\_print).
assign(max\_mem,64000).
assign(max\_weight,40).
assign(pick\_given\_ratio,6).

list(usable).
0 [] x=x.
end\_of\_list.

list(sos).
0 [] k x y=x.
0 [] abst x y z=x (k z) (y z).
0 [] y b(y) c(y)!=b(y) b(y) |\$ans(y).
end\_of\_list.

----> UNIT CONFLICT at   0.75 sec ----> 111 [binary,110.1,1.1] 

\$ans(abst abst k).

---------------- PROOF ----------------

1 [] x=x.
3,2 [] k x y=x.
4 [] abst x y z=x (k z) (y z).
5 [] x b(x) c(x)!=b(x) b(x) |\$ans(x).
6 [] x (k y) (z y)=abst x z y.
11 [para\_from,4.1.1,5.1.1.1]
 x (k b(abst x y)) (y b(abst x y)) c(abst x y)!=b(abst x y) b(abst x y) |
 \$ans(abst x y).
26 [para\_into,6.1.1.2,2.1.1] abst x (k y) z=x (k z) y.
110 [para\_from,26.1.1,11.1.1,demod,3]
 b(abst abst k) b(abst abst k)!=b(abst abst k) b(abst abst k) |\$ans(abst abst k).
111 [binary,110.1,1.1] \$ans(abst abst k).

------------ end of proof -------------
}
\tenpoint
\bigskip

Encouraged by this partial success we added the following axiom to 
the input file  
\roster
\item"{IX. \ }" $F\,x\,y=x\,x$
\endroster
and started the search for a contradiction anew.  While \O now
generated a number of properties of the combinator $F$ that to a
trained eye looked suspicious, it still did not find a contradiction.
We therefore decided for a different tack.  A contradiction should
have the form of Russell's paradox and should involve a contradictory
use of the equality combinator.  So we asked \O to find a combinator
$s$ with the property $s= Eq\l s,p_2\r$.  It is clear that such a 
self-referential combinator yields a contradiction: \ $s=p_1$ \iff
$s=p_2$.

Again, \O failed to find an answer, and again, we tried various
modifications, finally succeeding to elicit an answer when asked to
find an $s$ \st
$$
	s=Eq\l Ks, Kp_2\r.
$$
\O found such an $s$ in 29 sec (again, with options
{\tt para\_from\_units\_only} and {\tt para\_into\_units\_only}):  

\bigskip
\eightpoint

{\tt\obeylines\obeyspaces

set(knuth\_bendix).
set(ur\_res).
set(unit\_deletion).
set(para\_from\_units\_only).
set(para\_into\_units\_only).
set(bird\_print).
assign(max\_mem,64000).
assign(max\_weight,20).
assign(pick\_given\_ratio,6).

list(usable).
0 [] x=x.
end\_of\_list.

list(sos).
0 [] k x y=x.
0 [] p1 pair(x,y)=x.
0 [] p2 pair(x,y)=y.
0 [] pair(p1 x,p2 x)=x.
0 [] pair(x,y) z=pair(x z,y z).
0 [] abst x y z=x (k z) (y z).
0 [] eq pair(x,x)=p1.
0 [] x=y|eq pair(x,y) =p2.
0 [] x=y|x n(x,y) !=y n(x,y).
0 [] p1!=p2.
0 [] F x y=x x.
0 [] y!=eq pair(k y,k p2) |\$ans(y).
end\_of\_list.

----> UNIT CONFLICT at  28.79 sec ----> 3831 [binary,3830.1,3698.1] 
\$ans(abst (abst (abst (k eq))) pair(F,k (k p2)) 
(abst (abst (abst (k eq))) pair(F,k (k p2)))).

---------------- PROOF ----------------

3,2 [] k x y=x.
6 [] p2 pair(x,y)=y.
10 [] pair(x y,z y)=pair(x,z) y.
12 [] abst x y z=x (k z) (y z).
16 [] x=y|x n(x,y) !=y n(x,y).
18 [] F x y=x x.
20 [] eq pair(k x,k p2)!=x|\$ans(x).
23 [] x (k y) (z y)=abst x z y.
33 [para\_into,10.1.1.2,18.1.1] pair(x y,z z)=pair(x,F z) y.
37 [para\_into,10.1.1.2,6.1.1] pair(x pair(y,z),z)=pair(x,p2) pair(y,z).
42,41 [para\_into,10.1.1.2,2.1.1] pair(x y,z)=pair(x,k z) y.
50 [back\_demod,37,demod,42] pair(x,k y) pair(z,y)=pair(x,p2) pair(z,y).
54 [back\_demod,33,demod,42] pair(x,k (y y)) z=pair(x,F y) z.
59 [back\_demod,20,demod,42] eq (pair(k,k (k p2)) x)!=x|\$ans(x).
101 [para\_from,41.1.1,6.1.1.2] p2 (pair(x,k y) z)=y.
195 [para\_into,23.1.1.1,2.1.1] abst (k x) y z=x (y z).
206 [para\_into,23.1.1,12.1.1,demod,3] x (k (y z)) z=abst (abst x) y z.
1501,1500 [para\_from,50.1.1,101.1.1.2] p2 (pair(x,p2) pair(y,z))=z.
1884 [ur,54,16] pair(x,k (y y))=pair(x,F y).
1978,1977 [para\_from,1884.1.1,1500.1.1.2.2,demod,1501] k (x x)=F x.
2012 [para\_into,1977.1.1.2,18.1.1,demod,1978] F (F x)=F x.
2014 [para\_into,1977.1.1.2,2.1.1] F (k x)=k x.
2049,2048 [para\_from,2012.1.1,41.1.1.1,demod,42] pair(F,k x) (F y)=pair(F,k x) y.
2051,2050 [para\_from,2014.1.1,41.1.1.1,demod,42] pair(k,k x) y=pair(F,k x) (k y).
2074 [back\_demod,59,demod,2051] eq (pair(F,k (k p2)) (k x))!=x|\$ans(x).
3698 [para\_into,206.1.1,195.1.1,demod,3] abst (abst (abst (k x))) y z=x (y z).
3830 [para\_into,2074.1.1.2.2,1977.1.1,demod,2049]
 x x!=eq (pair(F,k (k p2)) x) |\$ans(x x).
3831 [binary,3830.1,3698.1] 
\$ans(abst (abst (abst (k eq))) pair(F,k (k p2)) 
(abst (abst (abst (k eq))) pair(F,k (k p2)))).

------------ end of proof -------------
}
\tenpoint
\bigskip

When inspecting
the proof we find that some of the axioms are not used.  After 
deleting the unnecessary axioms we obtained the following proof: \newline

\bigskip
\eightpoint

{\tt\obeylines\obeyspaces

set(knuth\_bendix).
set(ur\_res).
set(unit\_deletion).
set(para\_from\_units\_only).
set(para\_into\_units\_only).
set(bird\_print).
assign(max\_mem,64000).
assign(max\_weight,20).
assign(pick\_given\_ratio,6).

list(usable).
0 [] x=x.
end\_of\_list.

list(sos).
0 [] k x y=x.
0 [] p1 pair(x,y)=x.
0 [] p2 pair(x,y)=y.
0 [] pair(p1 x,p2 x)=x.
0 [] pair(x,y) z=pair(x z,y z).
0 [] abst x y z=x (k z) (y z).
0 [] x=y|x n(x,y) !=y n(x,y).
0 [] F x y=x x.
0 [] y!=eq pair(k y,k p2) |\$ans(y).
end\_of\_list.

----> UNIT CONFLICT at  34.19 sec ----> 4003 [binary,4002.1,158.1]

 \$ans(abst (k eq) pair(F,k (k p2)) (abst (k eq) pair(F,k (k p2)))).

---------------- PROOF ----------------

2 [] k x y=x.
7,6 [] p2 pair(x,y)=y.
10 [] pair(x y,z y)=pair(x,z) y.
13 [] x=y|x n(x,y) !=y n(x,y).
14 [] F x y=x x.
16 [] eq pair(k x,k p2)!=x|\$ans(x).
19 [] x (k y) (z y)=abst x z y.
27 [para\_into,10.1.1.2,14.1.1] pair(x y,z z)=pair(x,F z) y.
34,33 [para\_into,10.1.1.2,2.1.1] pair(x y,z)=pair(x,k z) y.
44 [back\_demod,27,demod,34] pair(x,k (y y)) z=pair(x,F y) z.
49 [back\_demod,16,demod,34] eq (pair(k,k (k p2)) x)!=x|\$ans(x).
158 [para\_into,19.1.1.1,2.1.1] abst (k x) y z=x (y z).
1468 [ur,44,13] pair(x,k (y y))=pair(x,F y).
1540,1539 [para\_from,1468.1.1,6.1.1.2,demod,7] k (x x)=F x.
1575 [para\_into,1539.1.1.2,14.1.1,demod,1540] F (F x)=F x.
1577 [para\_into,1539.1.1.2,2.1.1] F (k x)=k x.
1619,1618 [para\_from,1575.1.1,33.1.1.1,demod,34]
 pair(F,k x) (F y)=pair(F,k x) y.
1627,1626 [para\_from,1577.1.1,33.1.1.1,demod,34]
 pair(k,k x) y=pair(F,k x) (k y).
1661 [back\_demod,49,demod,1627] eq (pair(F,k (k p2)) (k x))!=x|\$ans(x).
4002 [para\_into,1661.1.1.2.2,1539.1.1,demod,1619]
 x x!=eq (pair(F,k (k p2)) x) |\$ans(x x).
4003 [binary,4002.1,158.1]
 \$ans(abst (k eq) pair(F,k (k p2)) (abst (k eq) pair(F,k (k p2)))).

------------ end of proof -------------

}
\tenpoint
\bigskip

This time, it took \O a little more time, but the resulting
combinator looks less complicated:  \newline
$$
	s=Abst(K\,Eq)\l F,K(K\,p_2)\r
	(Abst(K\,Eq)\l F,K(K\,p_2)\r).
$$
Upon closer inspection, the difference from the first proof is that
$Abst(K\,Eq)$ replaces $Abst(Abst(Abst(K\,Eq)))$.  This of course
begs the question whether these two terms might be equal.  Indeed,
the answer is yes, and in fact it is an instance of a general identity.
(The general identity $Abst(Abst(Abst\;x)))=Abst\;x$ is true in TRC
and is verified in the following section.)

Another observation is that $s=W\,W$ where $W$ is the combinator \linebreak
$Abst(K\, Eq)\l F,K(K\,p_2)\r$.  Thus the contradictory combinator
$W$ has the property
$$
	W\,W=Eq\l K(W\,W),K\,p_2\r,
$$
(which I believe should look familiar to experts in lambda calculus).

\medpagebreak

As a final experiment, we ask \O to verify that the combinator 
$s$ is contradictory.  So we add the axiom
\roster
\item"{X. \ }" $s=Eq \l K\,s,K\,p_2\r$
\endroster
and add a  name for the term $K\,p_2$.  We also disable the
{\tt para\_from\_units\_only} and {\tt para\_into\_units\_only} flags, as the more
general paramodulation is needed to produce a contradiction.

\bigskip
\eightpoint

{\tt\obeylines\obeyspaces

set(knuth\_bendix).
set(ur\_res).
set(unit\_deletion).
set(bird\_print).
assign(max\_mem,64000).
assign(max\_weight,40).
assign(pick\_given\_ratio,6).

list(usable).
0 [] x=x.
end\_of\_list.

list(sos).
0 [] k x y=x.
0 [] p1 pair(x,y)=x.
0 [] p2 pair(x,y)=y.
0 [] pair(p1 x,p2 x)=x.
0 [] pair(x,y) z=pair(x z,y z).
0 [] abst x y z=x (k z) (y z).
0 [] eq pair(x,x)=p1.
0 [] x=y|eq pair(x,y) =p2.
0 [] x=y|x n(x,y) !=y n(x,y).
0 [] p1!=p2.
0 [] k p1=cp1.
0 [] k p2=cp2.
0 [] s=eq pair(k s,k p2).
end\_of\_list.

----> UNIT CONFLICT at   0.31 sec ----> 84 [binary,82.1,17.1] \$F.

---------------- PROOF ----------------

2 [] k x y=x.
14,13 [] eq pair(x,x)=p1.
15 [] x=y|eq pair(x,y) =p2.
17 [] p2!=p1.
21,20 [] k p2=cp2.
22 [demod,21] eq pair(k s,cp2)=s.
30,29 [para\_from,20.1.1,2.1.1.1] cp2 x=p2.
35 [para\_into,15.2.1,22.1.1] k s=cp2|s=p2.
74,73 [para\_from,35.1.1,2.1.1.1,demod,30] s=p2.
82 [back\_demod,22,demod,74,21,14,74] p2=p1.
84 [binary,82.1,17.1] \$F.

------------ end of proof -------------

}

\tenpoint
\bigskip
\specialhead
Some Theorems of TRC
\endspecialhead

While running the experiments, we have observed a number of
interesting clauses generated by \O.  In addition to producing
equations true in TRC, some of the output led us to formulate, and
then verify (or disprove) various conjectures in TRC.  Below we give a
sample of some theorems of TRC that we found interesting (it remains
to be seen how important these facts are).  We view this as a modest
prelude to a systematic study of TRC, using automated reasoning
techniques. (We have obtained a large number of interesting theorems
of TRC that we intend to present in a future paper.) 
We hope that the information so obtained might contribute
to the eventual proof of inconsistency of NF (or to the construction of
a model).

In Proposition 1 below, $Id$ stands for the {\it identity} combinator,
$Id\;x=x$.  Note that $Id=\l p_1,p_2\r$.  Clearly, 1a is a consequence
of 1b. (The referee pointed out that the proof of 1b yields the stronger
statement 1c.)  The proofs of Proposition 1b and 2 are \O's and are
reprinted below with her permission.

\proclaim
{Proposition 1} \ \roster\runinitem"{(a)}" \ $Abst\;Abst\;Abst\;Abst\;
Abst\;Abst=Id$,  
\item"{(b)}" \ $Abst\;Abst\;Abst\;Abst = K(K(Id))$.
\item"{(c)}" \ $Abst\;Abst\;Abst\;Abst\;x\;y = Id.$
\endroster
\endproclaim

\proclaim
{Proposition 2} \ For all $x$, 
\roster
\item"{(a)}" \ $Abst(Abst(Abst\,x))=Abst\,x$,
\item"{(b)}" \ $Abst(Abst\,K(x))=K(x)$,
\item"{(c)}" \ $Abst\,K(K(x))=K(K(x))$.
\endroster
\endproclaim

[Warren Wood pointed out that 2c is also true in classical combinatory logic
(replacing $Abst$ with $S$ and $K(x)$ with $K$).]

\demo{Proofs of 1 (b), 2 (a), (b), (c)}

\bigskip
\eightpoint

{\tt\obeylines\obeyspaces

set(knuth\_bendix).
set(ur\_res).
set(unit\_deletion).
set(para\_from\_units\_only).
set(bird\_print).
assign(max\_mem,16000).
assign(max\_weight,40).
assign(pick\_given\_ratio,6).

list(usable).
0 [] x=x.
end\_of\_list.

list(sos).
0 [] k(x) y=x.
0 [] abst x y z=x k(z) (y z).
0 [] x=y|x n(x,y) !=y n(x,y).
0 [] id x=x.
0 [] abst abst abst abst!=k(k(id)).
end\_of\_list.

----> UNIT CONFLICT at   6.27 sec ----> 430 [binary,428.1,44.1] \$F.

---------------- PROOF ----------------

1 [] x=x.
3,2 [] k(x) y=x.
4 [] abst x y z=x k(z) (y z).
5 [] x=y|x n(x,y) !=y n(x,y).
6 [] id x=x.
8 [] k(k(id))!=abst abst abst abst.
9 [] x k(y) (z y)=abst x z y.
10 [ur,5,8,demod,3] abst abst abst abst n(k(k(id)),abst abst abst abst)!=k(id).
11 [para\_into,5.2.1,6.1.1] id=x|x n(id,x) !=n(id,x).
12 [para\_into,5.2.1,2.1.1] k(x)=y|y n(k(x),y) !=x.
17 [para\_from,4.1.1,11.2.1]
 abst x y=id|x k(n(id,abst x y)) (y n(id,abst x y)) !=n(id,abst x y).
21 [para\_into,9.1.1.1,6.1.1,demod,3] abst id x y=y.
38,37 [para\_from,9.1.1,4.1.1.1,demod,3] abst abst x y z=y (x y z).
39 [back\_demod,10,demod,38] abst (abst abst n(k(k(id)),abst abst abst abst))!=k(id).
40 [ur,21,11] abst id x=id.
43,42 [ur,40,12] k(id)=abst id.
44 [back\_demod,39,demod,43,43]
 abst (abst abst n(k(abst id),abst abst abst abst))!=abst id.
51 [para\_from,40.1.1,5.2.1] abst id=x|x n(abst id,x) !=id.
407 [para\_from,37.1.1,17.2.1,demod,3,unit\_del,1] abst (abst abst x) y=id.
428 [ur,407,51] abst (abst abst x)=abst id.
430 [binary,428.1,44.1] \$F.

------------ end of proof -------------
}

\bigskip
{\tt\obeylines\obeyspaces

set(knuth\_bendix).
set(ur\_res).
set(unit\_deletion).
set(para\_from\_units\_only).
set(bird\_print).
assign(max\_mem,64000).
assign(max\_weight,60).
assign(pick\_given\_ratio,6).

list(usable).
0 [] x=x.
end\_of\_list.

list(sos).
0 [] k(x) y=x.
0 [] abst x y z=x k(z) (y z).
0 [] x=y|x n(x,y) !=y n(x,y).
0 [] abst (abst (abst b))!=abst b.
end\_of\_list.

----> UNIT CONFLICT at   4.66 sec ----> 638 [binary,636.1,8.1] \$F.

---------------- PROOF ----------------

3,2 [] k(x) y=x.
4 [] abst x y z=x k(z) (y z).
5 [] x=y|x n(x,y) !=y n(x,y).
6 [] abst (abst (abst b))!=abst b.
7 [] x k(y) (z y)=abst x z y.
8 [ur,5,6] abst (abst (abst b)) n(abst (abst (abst b)),abst b)
!=abst b n(abst (abst (abst b)),abst b).
20 [para\_into,7.1.1.2,2.1.1] x k(y) z=abst x k(z) y.
21 [para\_into,7.1.1,4.1.1,demod,3] x k(y z) z=abst (abst x) y z.
48 [para\_into,20.1.1,7.1.1] abst x k(y z) z=abst x y z.
457 [para\_into,48.1.1,21.1.1] abst (abst (abst x)) y z=abst x y z.
636 [ur,457,5] abst (abst (abst x)) y=abst x y.
638 [binary,636.1,8.1] \$F.

------------ end of proof -------------
}

\bigskip
{\tt\obeylines\obeyspaces

set(knuth\_bendix).
set(ur\_res).
set(unit\_deletion).
set(para\_from\_units\_only).
set(bird\_print).
assign(max\_mem,16000).
assign(max\_weight,40).
assign(pick\_given\_ratio,6).

list(usable).
0 [] x=x.
end\_of\_list.

list(sos).
0 [] k(x) y=x.
0 [] abst x y z=x k(z) (y z).
0 [] x=y|x n(x,y) !=y n(x,y).
0 [] id x=x.
0 [] k(b)!=abst (abst k(b)).
end\_of\_list.

----> UNIT CONFLICT at   1.11 sec ----> 104 [binary,102.1,10.1] \$F.

---------------- PROOF ----------------

3,2 [] k(x) y=x.
5 [] x=y|x n(x,y) !=y n(x,y).
8 [] abst (abst k(b))!=k(b).
9 [] x k(y) (z y)=abst x z y.
10 [ur,5,8,demod,3] abst (abst k(b)) n(abst (abst k(b)),k(b))!=b.
24 [para\_into,9.1.1.1,2.1.1] abst k(x) y z=x (y z).
68 [para\_into,24.1.1,9.1.1,demod,3] abst (abst k(x)) y z=x z.
102 [ur,68,5] abst (abst k(x)) y=x.
104 [binary,102.1,10.1] \$F.

------------ end of proof -------------

}

\bigskip
{\tt\obeylines\obeyspaces

set(knuth\_bendix).
set(ur\_res).
set(unit\_deletion).
set(para\_from\_units\_only).
set(bird\_print).
assign(max\_mem,16000).
assign(max\_weight,40).
assign(pick\_given\_ratio,6).

list(usable).
0 [] x=x.
end\_of\_list.

list(sos).
0 [] k(x) y=x.
0 [] abst x y z=x k(z) (y z).
0 [] x=y|x n(x,y) !=y n(x,y).
0 [] id x=x.
0 [] k(k(b))!=abst k(k(b)).
end\_of\_list.

----> UNIT CONFLICT at   0.53 sec ----> 55 [binary,54.1,1.1] \$F.

---------------- PROOF ----------------

1 [] x=x.
3,2 [] k(x) y=x.
5 [] x=y|x n(x,y) !=y n(x,y).
8 [] abst k(k(b))!=k(k(b)).
9 [] x k(y) (z y)=abst x z y.
10 [ur,5,8,demod,3] abst k(k(b)) n(abst k(k(b)),k(k(b)))!=k(b).
25,24 [para\_into,9.1.1.1,2.1.1] abst k(x) y z=x (y z).
54 [ur,10,5,demod,25,3,3] b!=b.
55 [binary,54.1,1.1] \$F.

------------ end of proof -------------

}
\enddemo
\tenpoint

\enddocument